\documentclass[graybox]{svmult}

\usepackage{mathptmx}       
\usepackage{helvet}         
\usepackage{courier}        
\usepackage{type1cm}        
\usepackage{makeidx}         
\usepackage{graphicx}        
\usepackage{multicol}        
\usepackage[bottom]{footmisc}

\usepackage{amsmath,amssymb}

\usepackage{tikz}

\makeindex             

\newtheorem{assumption}{Assumption}[section]
\colorlet{mygreen}{green!50!gray}
\DeclareMathOperator{\e}{e}

\begin{document}

\title*{Asymptotic-preserving methods for an anisotropic model of electrical potential in a tokamak}
\titlerunning{AP methods for an anisotropic model of electrical potential in a tokamak}
\author{Philippe Angot, Thomas Auphan and Olivier Gu\`es}
\institute{Philippe Angot, Thomas Auphan \and Olivier Gu\`es \at Aix Marseille Universit\'e, CNRS, Centrale Marseille, I2M, UMR 7373, 13453 Marseille, France \email{[thomas.auphan, philippe.angot, olivier.gues]@univ-amu.fr}}
%
%
\maketitle


\abstract{
A 2D nonlinear model for the electrical potential in the edge plasma in a tokamak generates a stiff problem due to the low resistivity in the direction parallel to the magnetic field lines.
An asymptotic-preserving method based on a micro-macro decomposition is studied in order to have a well-posed problem, even when the parallel resistivity goes to $0$. Numerical tests with a finite difference scheme show a bounded condition number for the linearised discrete problem solved at each time step, which confirms the theoretical analysis on the continuous problem.
}

\section{Introduction}
\label{sec:1}

The fusion reaction can be performed using a tokamak, a machine whose shape is toroidal. The plasma is confined and warmed in the core of the tokamak to produce the fusion reaction. This technique is expected to maintain the fusion reaction during a long time (more than five minutes, for the ITER project).

One of the main challenges for this objective is to control the wall-plasma interactions. Indeed, the magnetic confinement is not perfect and the plasma is in contact with the wall. In a tokamak such as TORE SUPRA, an obstacle called the limiter, is settled at the bottom of the machine. Due to the strong magnetic confinement, the plasma transport essentially occurs along the magnetic field lines. Thus, the parallel resistivity $\eta$ is  very small (typically, $\eta =10^{-6}$), generating a strong anisotropy in the model. The area where the magnetic lines are interrupted by the limiter is called the scrape-off layer. The numerical simulation of the edge plasma transport allows us to better understand the interactions with the wall.

\section{Anisotropic model of the electrical potential}
\begin{figure}
\begin{center}
\begin{tikzpicture}[scale=1.2]
\fill [color=gray!15] (0,0) rectangle (0.8,	1.5);
\fill [color=gray!15] (3.2,0) rectangle (4,1.5);
\draw (0,0) rectangle (4,2.5);
\draw (0.4,0) node[right,rotate=90]{Limiter};
\draw (3.6,0) node[right,rotate=90]{Limiter};
\draw [color=gray!90] (-0.4,0.1) node[right,rotate=90]{Periodic BC};
\draw [color=gray!90] (4.4,0.1) node[right,rotate=90]{Periodic BC};
\draw [thick, ->] (0,0) -- (4.5,0);
\draw [thick, ->] (0,0) -- (0,2.75);
\draw (0,0) node[below]{$-0.5$};
\draw (0,0) node[left]{$0$};
\draw (0,2.5) node[left]{$1$};
\draw (0.8,0) node[below]{$-L$};
\draw (3.2,0) node[below]{$L$};
\draw (4,0) node[below]{$0.5$};
\draw [thick, gray!70] (0.8,0)--(3.2,0);
\draw (2,0) node[above,gray!70]{$\Sigma_\parallel$};
\draw [thick, gray!70] (0,2.5)--(4,2.5);
\draw (2,2.5) node[below,gray!70]{$\Sigma_\parallel$};
\draw [very thick] (0.8,0)--(0.8,1.5);
\draw [very thick] (3.2,0)--(3.2,1.5);
\draw [thick, gray!70] (0,1.5)--(0.8,1.5);
\draw [thick, gray!70] (3.2,1.5)--(4,1.5);
\draw (0,1.5) node[left]{$l$};
\draw (3.6,1.5) node[above,gray!70]{$\Sigma_\parallel$};
\draw (0.4,1.5) node[above,gray!70]{$\Sigma_\parallel$};
\draw (4.6,0) node[below]{$x$};
\draw (0,2.7) node[left]{$y$};
\draw (2,1.25) node{Plasma};
\draw (2,0.95) node{$\Omega$};
\draw (2,0) node[below]{Wall};
\draw (2,2.5) node[above]{Center};
\end{tikzpicture}
\end{center}
\caption{Schematic representation of the 2D domain.}\label{Fig_domain}
\end{figure}
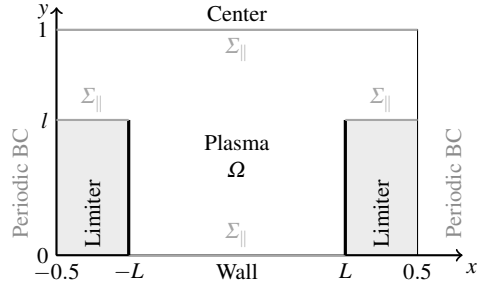

In this paper, we focus on a 2D model of the electrical potential of the edge plasma $\phi_\eta$ in a tokamak with a limiter configuration. A schematic representation of the domain is given in Fig. \ref{Fig_domain}. The $x$ axis corresponds to the curvilinear coordinates along a magnetic field line and the $y$ axis is the radial direction. In the following equations, the curvature terms have been neglected. As the magnetic field lines above the limiter set are closed, periodic boundary conditions are imposed at $x = \pm 0.5$.
%
The dimensionless problem for the electrical potential reads:
\begin{equation}\label{phi_2D_model}
\left\{\begin{aligned}
& -\partial_t \partial_y^2 \phi_\eta - \frac{1}{\eta} \partial_x^2 \phi_\eta + \nu \partial_y^4 \phi_\eta = S & &\text{in } ]0,T[ \times \Omega\\
&\partial_y \phi_{\eta|t=0} = \partial_y \phi_{ini} & &\text{in } \Omega\\
&\partial_y \phi_{\eta|\Sigma_\parallel} = 0 \quad \text{ and } \quad \partial^3_y \phi_{\eta|\Sigma_\parallel} = 0 & &\text{on } ]0,T[ \times \Sigma_\parallel\\
&\partial_x \phi_{\eta|x=-L} = \eta \left(1 - \e^{\Lambda - \phi_{\eta|x=-L}}\right) & &\text{on } ]0,T[ \times ]0,l[ \times \{ -L \}\\
&\partial_x \phi_{\eta|x=L} = - \eta \left(1 - \e^{\Lambda - \phi_{\eta|x=L}}\right) & & \text{on } ]0,T[ \times ]0,l[ \times \{ L \},
\end{aligned} \right.
\end{equation}
where $\nu$ corresponds to the ionic viscosity in the perpendicular direction and $\Lambda$ stands for the reference potential inside the limiter. The initial condition is $\partial_y \phi_{\eta|t=0} = \partial_y \phi_{ini}$.
Negulescu \emph{et al.} \cite{Neg08} proved that, for a fixed value of $\eta>0$, the problem (\ref{phi_2D_model}) admits a unique weak solution, under suitable hypotheses on the data $\phi_{ini}$ and $S$.

The boundary conditions at the limiter interface $x = \pm L$, are nonlinear. Setting directly $\eta=0$ in the system (\ref{phi_2D_model}) (after multiplying the first equation by $\eta$) leads to an under-determined problem since there are only homogeneous Neumann boundary conditions at the limiter surface $x=\pm L$. Thus, when $\eta$ is small the numerical resolution of the problem (\ref{phi_2D_model}) becomes stiff. This issue can be avoided by reformulating the problem (\ref{phi_2D_model}) thanks to asymptotic-preserving methods.

\section{The micro macro asymptotic-preserving method}
We study the Asymptotic-Preserving (AP) method introduced by Degond \emph{et al.} \cite{Deg12} for a linear anisotropic elliptic problem. It consists in a decomposition of the solution $\phi_\eta$ as $\phi_\eta = p_\eta + \eta q_\eta$ where $\partial_x p_\eta = 0$ and $q_{\eta | x = -L} = 0$. Then, it yields the problem below where the unknowns are $(\phi_\eta,q_\eta)$:
\begin{equation}\label{phi_2D_model_q}
\left\{\begin{aligned}
& -\partial_t \partial_y^2 \phi_\eta -  \partial_x^2 q_\eta + \nu \partial_y^4 \phi_\eta = S & &\text{in } ]0,T[ \times \Omega\\
& \partial_x^2 \phi_\eta = \eta \partial_x^2 q_\eta & &\text{in } ]0,T[ \times \Omega\\
&\partial_x \phi_{\eta|x=-L}=\eta \partial_x q_{\eta|x=-L} & &\text{on } ]0,T[ \times ]0,l[ \times \{ -L \}\\
&\partial_x \phi_{\eta|x=L}=\eta \partial_x q_{\eta|x=L} & & \text{on } ]0,T[ \times ]0,l[ \times \{ L \}\\
&\partial_x \phi_{\eta|x=-0.5}=\eta \partial_x q_{\eta|x=-0.5} & &\text{on } ]0,T[ \times ]l,1[ \times \{ -0.5 \}\\
&\partial_x \phi_{\eta|x=0.5}=\eta \partial_x q_{\eta|x=0.5} & & \text{on } ]0,T[ \times ]l,1[ \times \{ 0.5 \}\\
&\partial_y \phi_{\eta|t=0} = \partial_y \phi_{ini} & &\text{in } \Omega\\
&\partial_y \phi_{\eta|\Sigma_\parallel} = 0 \quad \text{ and } \quad \partial^3_y \phi_{|\Sigma_\parallel} = 0 & &\text{on } ]0,T[ \times \Sigma_\parallel\\
&\partial_x q_{\eta|x=-L} = \left(1 - \e^{\Lambda - \phi_{\eta|x=-L}}\right) & &\text{on } ]0,T[ \times ]0,l[ \times \{ -L \}\\
&\partial_x q_{\eta|x=L} = -\left(1 - \e^{\Lambda - \phi_{\eta|x=L}}\right) & & \text{on } ]0,T[ \times ]0,l[ \times \{ L \},
\end{aligned} \right.
\end{equation}
One important advantage of this AP method is that it can be easily implemented even if the mesh is not aligned with the directions $(O x)$ and $(O y)$. The main drawback is the need to compute two unknowns ($\phi_\eta$ and $q_\eta$) on the 2D domain though only $\phi_\eta$ is interesting for the physics.

Let us give the theoretical result which ensures that the modified problem is well-posed for $\eta=0$, and that $\phi_\eta$ converges towards $\phi_0$. First, we provide the definitions of the spaces used for the variational formulation of the problem (\ref{phi_2D_model_q}).

\begin{definition}
Let us define the following Hilbert spaces:
\begin{itemize}
 \item $V=\left\{f \in H^1(\Omega), \partial_y^2 f \in L^2(\Omega), f \text{ periodic on } \{-0.5,0.5\} \times ]l,1[, \partial_y f = 0 \text{ on } {\Sigma_\parallel} \right\}$ with the scalar product:
 \begin{equation*}
\langle f,u \rangle_{V} = \int_\Omega{\partial_x f \, \partial_x u \, dy dx}+\int_\Omega{\partial^2_y f \, \partial^2_y u \, dy dx} + 2 \int_0^l{f_{|x=L} \, u_{|x=L} \, dy}.
\end{equation*}
 \item $Q = \left\{f \in L^2(\Omega), \partial_x f \in L^2(\Omega) , f_{|x=-L}=0 \text{ on } ]0,1[\right\}$, with the scalar product:
\begin{equation*} 
 \langle f,u\rangle_Q = \int_\Omega{\partial_x f \, \partial_x u \, dy dx}.
 \end{equation*}
\end{itemize}
\end{definition}


\begin{definition}
The space $\mathcal{A}$ is the set of functions $\phi$ such that:
\begin{itemize}
 \item $\phi \in L^2(0,T;V)$.
 \item $\partial_y \phi \in L^\infty(0,T;L^2(\Omega))$.
 \item $\partial_y \phi \in L^2\left(0,T;\{f \in H^1(\Omega), \partial_y^2 f \in L^2(\Omega), f_{|\Sigma_\parallel}=0\}\right)$.
 \item $\partial_y^2 \phi \in L^\infty(0,T;L^2(\Omega))$.
 \item $\partial_t \phi \in L^2(0,T;V)$.
 \item $\partial_y \partial_t \phi \in L^\infty(0,T;L^2(\Omega))$.
\end{itemize}
\end{definition}
The weak solution $\phi_\eta$ of (\ref{phi_2D_model_q}) is then searched in the space $\mathcal{A}$.

\begin{assumption}\label{Intro_hypothese_phi_2D}
Assume that $S$ and $\phi_{ini}$ verify:
\begin{enumerate}
 \item $S, \partial_y S, \partial_y^2 S, \partial_t S, \partial_t^2 S \in L^2(]0,T[ \times \Omega)$, $\| S \|_{L^\infty(]0,T[ \times \Omega)} \leq C_s$ and $\| S_{|t=T} \|_{L^\infty(\Omega)} \leq C_s$ with $C_s$ sufficiently small.
 \item $\phi_{ini} \in H^4(\Omega)$.
 \item $\phi_{ini}$ does not depend on $x$.
 \item  $\displaystyle\int_\Omega{S_{|t=0}\, dy dx} = \nu \displaystyle\int_\Omega{\partial_y^4 \phi_{ini} \, dy dx} + 2 \displaystyle\int_0^l{\left( 1 -  \e^{\Lambda - \phi_{ini|x=L}}\right) \, dy}$.
\end{enumerate}
\end{assumption}
The two last hypotheses are compatibility conditions for the initial and boundary conditions with the source term.

We can now write the theorem which asserts the convergence of $\phi_\eta$ to $\phi_0$ when $\eta$ goes to $0$:
\begin{theorem}\label{Intro_phi_2D_th}
With the assumption \ref{Intro_hypothese_phi_2D}, the weak formulation of  (\ref{phi_2D_model_q}):\\ 
find $(\phi_\eta,q_\eta) \in \mathcal{A} \times L^2(0,T;Q)$ verifying
\begin{equation}\label{Intro_phi_2D_AP_pb_non_lin_faible}
\left\{ \begin{aligned}
&\forall \xi \in H^1(]0,T[), \forall u \in V \cap H^2(\Omega), \forall w \in Q,\\
& \int_\Omega{\partial_y \phi_{\eta|t=T} \, \partial_y u \, dy dx} \xi(T) - \int_0^T{\int_\Omega{ \partial_y \phi_{\eta|t=T} \,  \partial_y u \, dy dx} \, \xi' \, dt} \\ 
& \quad +  \int_0^T{\int_\Omega{\partial_x q_\eta \, \partial_x u\,  dy dx}\,  \xi dt}  + \nu \int_0^T{\int_\Omega{\partial_y^2 \phi_\eta \, \partial_y^2 u \, dy dx}\, \xi \, dt} \\
& \quad + \int_0^T{\int_0^l{ \left(1-\e^{\Lambda-\phi_{\eta|x=-L}} \right) \, u_{|x=-L} \, dy}\, \xi \, dt} + \int_0^T{\int_0^l{ \left(1-\e^{\Lambda-\phi_{\eta|x=L}} \right) \, u_{|x=L} \, dy}\, \xi \, dt} \\
&\quad =\int_\Omega{ \partial_y \phi_{ini} \, \partial_y u \, dy dx}\, \xi(0) + \int_0^T{\int_\Omega{S \, u \, dy dx}\, \xi \, dt}\\
& \eta \int_0^T{\int_\Omega{\partial_x q_\eta \, \partial_x w \, dy dx} \, \xi \, dt} =\int_0^T{\int_\Omega{\partial_x \phi_\eta \, \partial_x w \, dy dx} \, \xi \, dt},
 \end{aligned}\right.
\end{equation}
admits a unique solution. Besides, $(\phi_\eta, q_\eta)$ converges weakly in $L^2(]0,T[\times\Omega)^2$, towards $(\phi_0,q_0) \in \mathcal{A} \times L^2(0,T;Q)$ the solution of (\ref{Intro_phi_2D_AP_pb_non_lin_faible}) when $\eta$ equals $0$.

Finally, the following error estimate holds:
\begin{equation*}
 \|\phi_\eta - \phi_0\|_{L^1(0,T;L^2(\Omega))} \leq c(T,\Omega,\phi_0,S,\Lambda) \, \sqrt{\eta},
\end{equation*}
where $c(T,\Omega,\phi_0,S,\Lambda)>0$ does not depend on $\eta$.
\end{theorem}
Theorem \ref{Intro_phi_2D_th} provides an error estimate for the norm in $L^1(0,T;L^2(\Omega))$, but not for the $L^2(]0,T[\times\Omega)$ norm. This point can be subject to further improvements.

This result is shown in \cite{Ang14_2,Aup14}. The proof of the existence and uniqueness of $\phi_0$ follows the same steps of \cite{Neg08}, based on a fixed point method. The existence and uniqueness of $q_0$ and the convergence of $(\phi_\eta,q_\eta)$ when $\eta$ goes to $0$ are shown by extending to a nonlinear case the proof provided in \cite{Deg12} for a linear elliptic problem.

\section{Numerical experiments}
In this section, some numerical tests are presented for the system (\ref{phi_2D_model_q}).
The space discretisation is done by the centred finite difference scheme. The time resolution uses Euler semi-implicit method.

At first glance, a directional splitting method seems to be interesting. But, the discrete problems obtained in the directions $x$ and $y$ are not invertible.
The problem is thus discretised implicitly, except for the nonlinear term.
At each time step, a linear system has to be solved to compute the approximations of $\phi_\eta$ and $q_\eta$. 


Let us consider a rectangular mesh of the space domain $\Omega$ with a constant mesh step $\delta x$ (for the direction $(O x)$) and $\delta y$ (for the direction $(O y)$). The time step writes $\delta t$.
The scalar quantities $\phi_{i,j}^n, q_{i,j}^n$ stands respectively for the approximations of $\phi_\eta(n \delta t, -0.5+i \delta x, j \delta y)$ and  
$q_\eta(n \delta t, -0.5 + i \delta x, j \delta y)$. The boundary condition at $x=-L$ is discretised as:
\begin{equation*}
\frac{q_{I_1+1,j}^{n+1} - q_{I_1-1,j}^{n+1}}{2 \delta x} - \phi_{I_1,j}^{n+1}= \left(1- \e^{\Lambda- \phi_{I_1,j}^n} - \phi_{I_1,j}^n\right),
\end{equation*}
where $I_1$ is the index such that $-0.5+ I_1 \delta x = -L$.

For the boundary condition at $x=L$, the same technique is used.
This time linearisation enables us to have an invertible matrix which is the same at each time step.

The mesh convergence test is performed using a configuration where the limiter goes up to the top of the computational domain, \emph{i.e.} $l=1$. This does not change the results proven for $l<1$.  For $L=0.4$, the chosen manufacturated solution is 
\begin{equation}\label{Phi_2D_AP_sol_test}
 \phi_\eta(t,x,y)=\eta \left(\dfrac{t}{\pi}\right)^2 \cos(\pi y) \cos(1.25 \pi x) - \ln\left(1-\dfrac{1.25 t^2}{\pi \cos(\pi y)}\right) + \Lambda.
\end{equation}
Let us note that the source term $S$ associated to the manufactured solution (\ref{Phi_2D_AP_sol_test}) depends on $\eta$ but is not singular when $\eta$ goes to $0$. This differs from the hypotheses made for Theorem \ref{Intro_phi_2D_th}.

\begin{figure}
\begin{center}
\includegraphics[scale=0.3, trim = 4mm 60mm 15mm 85mm, clip=true] {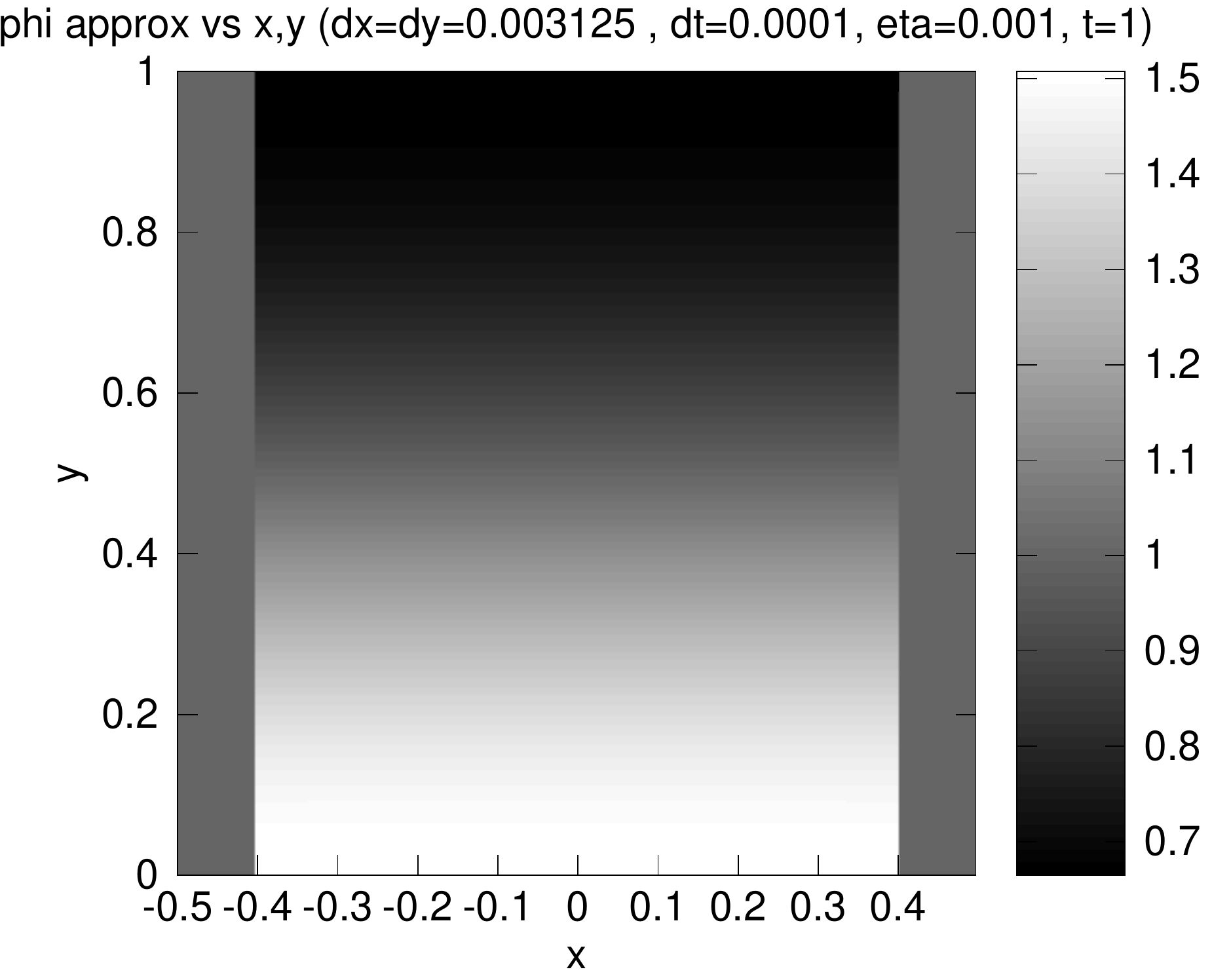}
\includegraphics[scale=0.3, trim = 6mm 60mm 14mm 85mm, clip=true]{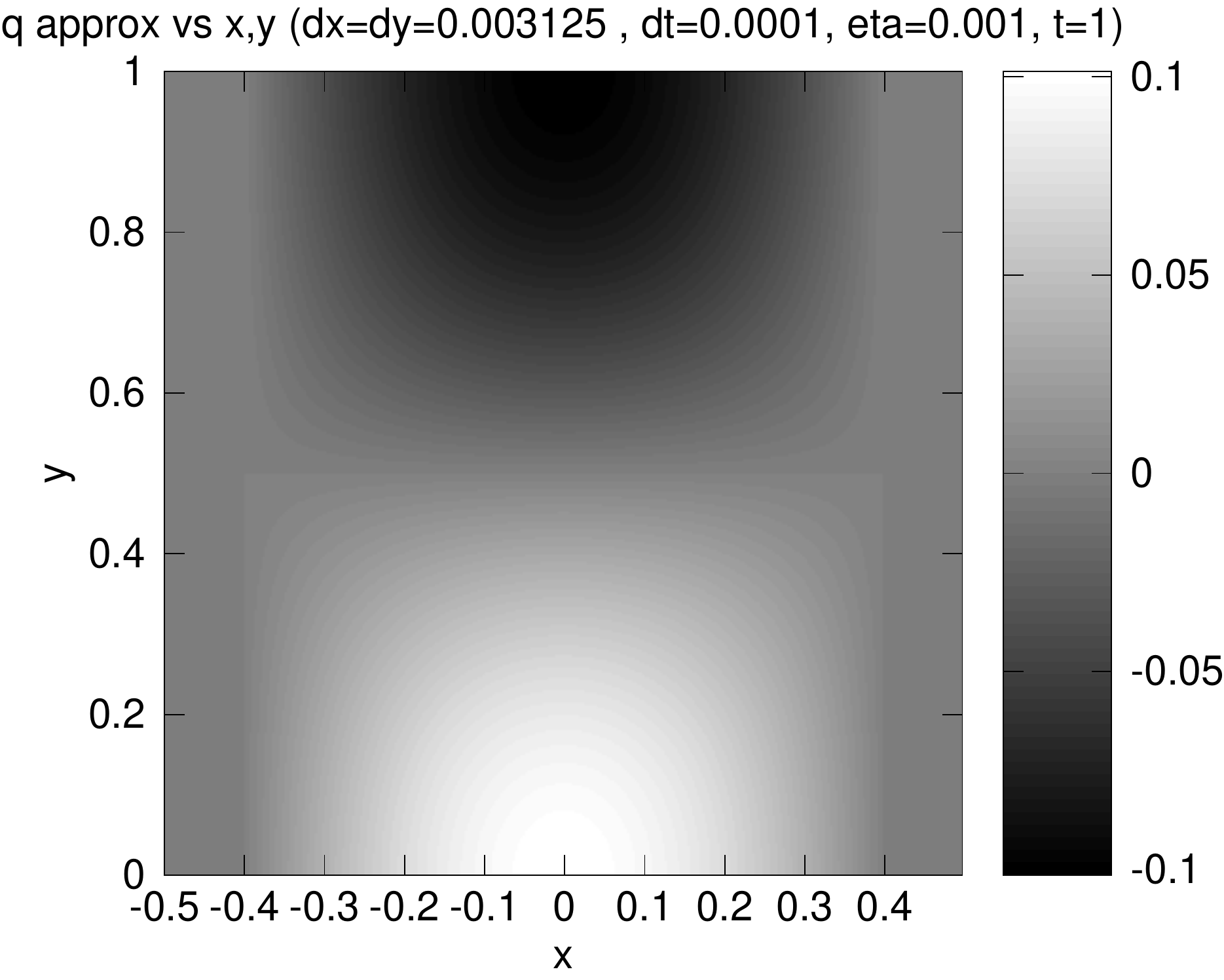}
\caption{Approximate fields of $\phi_\eta$ and $q_\eta$ for $\delta x = \delta y = 0.003125$, $\delta t=0.0001$ and $\eta=0.001$. The reference solution is given by (\ref{Phi_2D_AP_sol_test}). Recall that the limiter area corresponds to $x\leq -0.4$ and $x\geq 0.4$: the values of  $\phi_\eta$ do not have any physical sense in this zone.}\label{Cour_AP_2D_trace_fig} 
\end{center}
\end{figure}

\begin{figure}
\begin{center}
\includegraphics[scale=0.55, trim = 4mm 2mm 9mm 8mm, clip=true]{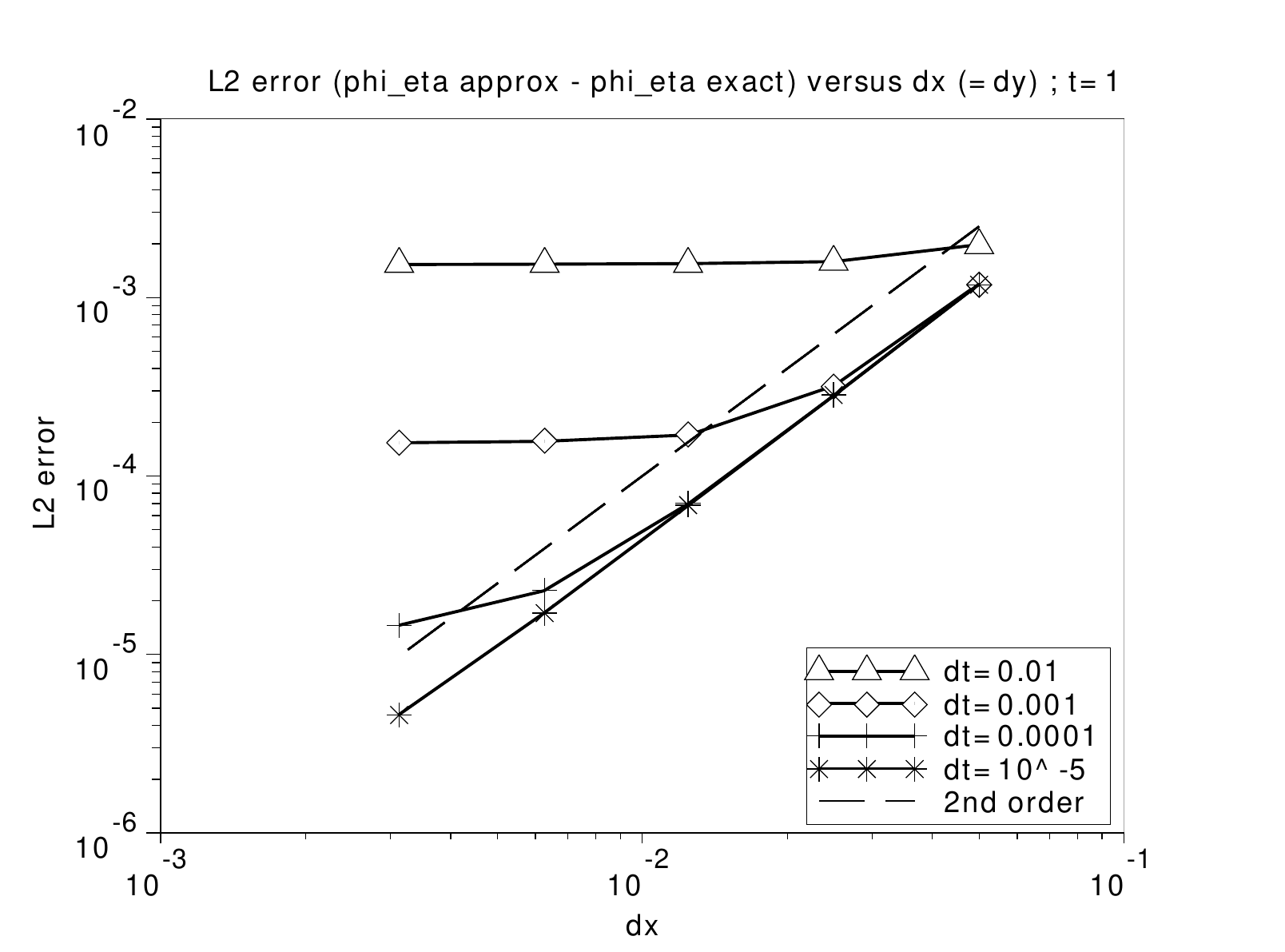}
\caption{$\|\phi_{\eta}^{approx} - \phi_\eta\|_{L^2(\Omega)}$ at $t=1$ as a function of the space step $\delta x = \delta y$  for different values of the time step and $\eta=0.001$. The reference solution is given by (\ref{Phi_2D_AP_sol_test}).}\label{Cour_AP_2D_err_L2_vs_dx_fig} 
\end{center}
\end{figure}

The plot of the approximated solution is shown in Fig. \ref{Cour_AP_2D_trace_fig}.
Studying the $L^2$ error in Fig. \ref{Cour_AP_2D_err_L2_vs_dx_fig}, we observe that the numerical scheme is of second-order accuracy in space. 

\begin{figure}
\begin{center}
\includegraphics[scale=0.55, trim = 4mm 2mm 9mm 8mm, clip=true]{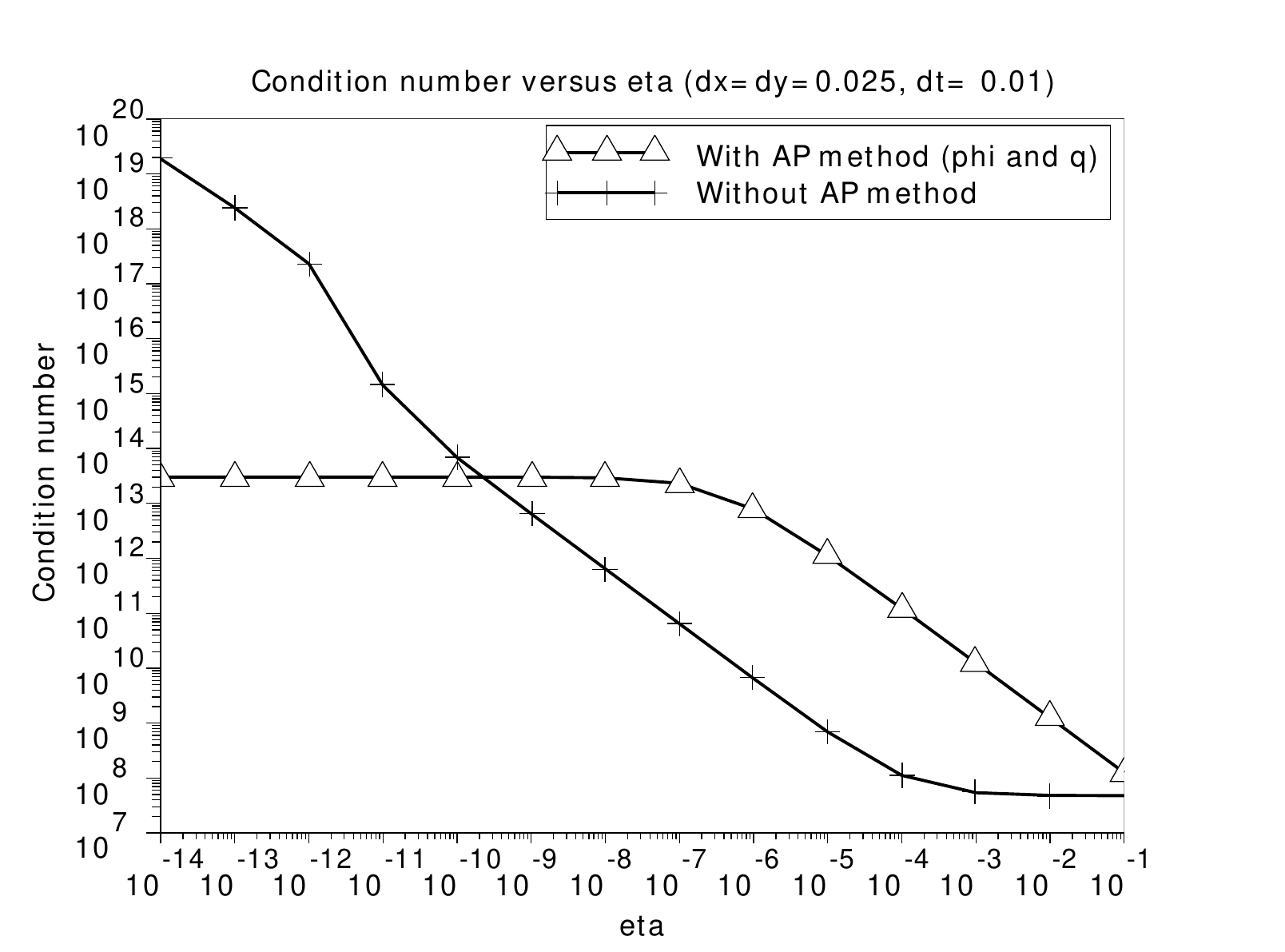}
\caption{Condition number in the Euclidean norm as a function of the parallel resistivity $\eta$ for the linear system approaching (the same at each time step) the solution (\ref{Phi_2D_AP_sol_test}) with $\delta x = \delta y = 0.025$ and $\delta t = 0.001$.
}\label{Cour_AP_2D_cond_fig} 
\end{center}
\end{figure}

In Fig. \ref{Cour_AP_2D_cond_fig}, we observe that the condition number obtained with the AP method is high but it is bounded independently from $\eta$. This is not the case for the matrix obtained for the resolution of (\ref{phi_2D_model}) without the asymptotic-preserving method. In order to avoid the issues due to the bad conditioning, we choose a LU method to solve the linear problem at each time step, which is faster than a GMRES solver with PETSc library. Finding an efficient preconditioner in order to use iterative methods is a future enhancement of this work.

\begin{figure}
\begin{center}
\includegraphics[scale=0.55, trim = 4mm 2mm 9mm 8mm, clip=true]{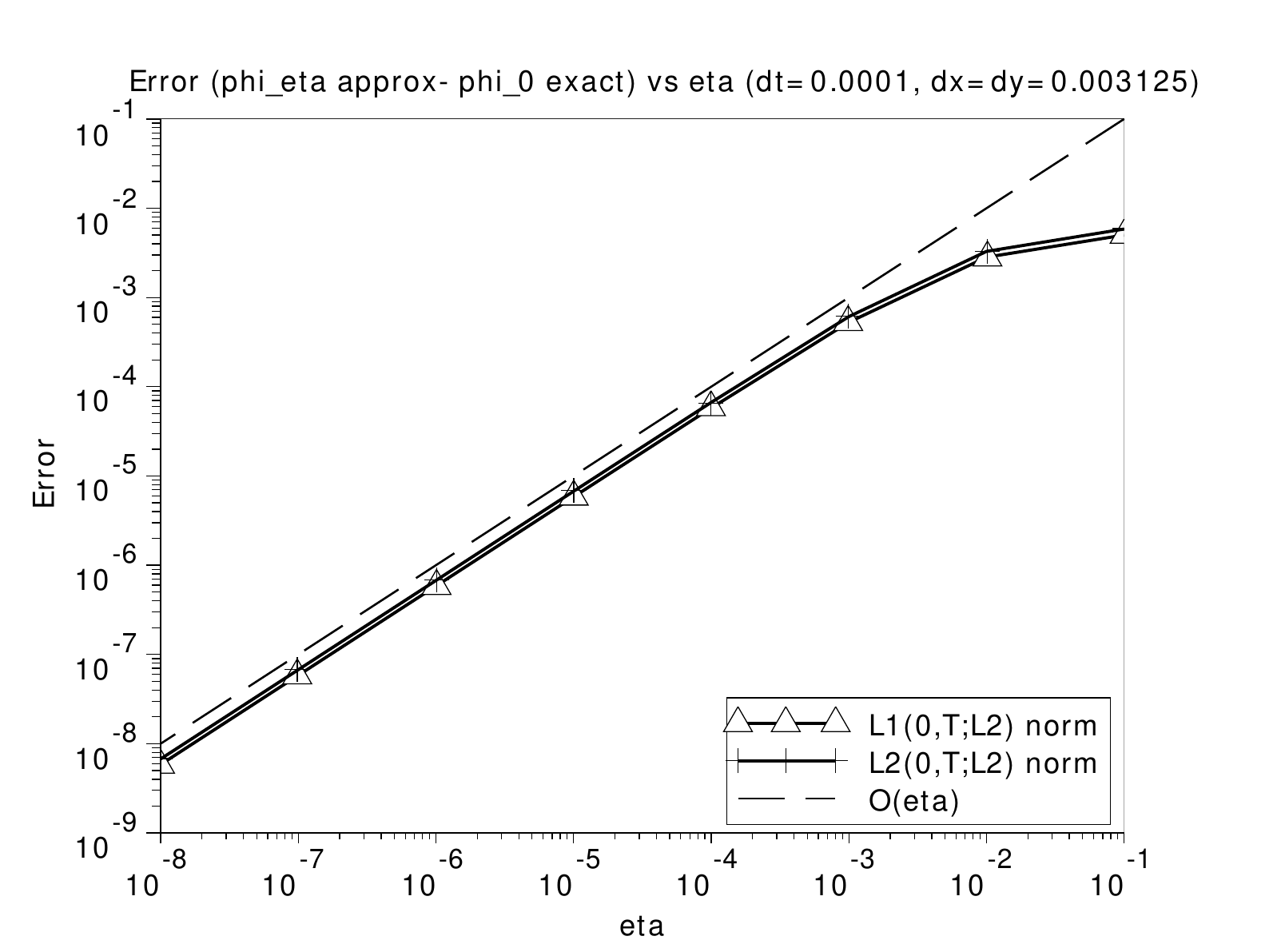}
\caption{$\|\phi_\eta - \phi_0\|_{L^1(0,T;L^2(\Omega))}$ ($\Delta$) and $\|\phi_\eta - \phi_0\|_{L^2(]0,T[ \times \Omega)}$ ($+$) as a function of $\eta$. The configuration is given by Eqs. (\ref{Phi_2D_AP_sol_test_2}) with $T=1$, $\delta x = \delta y = 0.003125$ and $\delta t=0.0001$.
}\label{Cour_AP_2D_cnv_eta_to_0_fig} 
\end{center}
\end{figure}

For the convergence when $\eta$ tends to $0$, the same domain is considered ($l=1, L=0.4$) but another source term is chosen :
\begin{equation}\label{Phi_2D_AP_sol_test_2}
S(t,x,y)= 40 \, t \cos(2\pi \, y) \, \sin\left(\dfrac{\pi}{2 L} \, x \right) \quad , \quad \phi_{ini}(x,y)=\Lambda=0
\end{equation}
This configuration (\ref{Phi_2D_AP_sol_test_2}) with $l=1$ leads to $\phi_0(t,x,y)=0$, which enables us to compute numerically $\|\phi_\eta - \phi_0\|_{L^1(0,T;L^2(\Omega))}$ and $\|\phi_\eta - \phi_0\|_{L^2(0,T;L^2(\Omega))}$. For these two norms, we observe a convergence in $\mathcal{O}(\eta)$, see Fig. \ref{Cour_AP_2D_cnv_eta_to_0_fig}. This suggests that the estimate of Theorem \ref{Intro_phi_2D_th} might be improved.


\section{Conclusion}
The high anisotropy of the 2D model for the edge plasma electrical potential in a tokamak leads to an ill-conditioned matrix for the numerical approximation using classical methods. The micro-macro decomposition induced by Degond \emph{et al.} \cite{Deg12} for a linear anisotropic elliptic problem is studied and analysed for the nonlinear evolution problem of the electrical potential. This method yields a weak formulation which is not degenerated when the parallel resistivity $\eta$ tends to $0$. Moreover, we have the estimate 
\begin{equation*}
\|\phi_\eta - \phi_0\|_{L^1(0,T,L^2(\Omega))} = \mathcal{O}\left(\sqrt{\eta}\right),
\end{equation*}
which can probably be improved, as suggested by the numerical results.

\noindent{\bf Acknowledgements:} This work has been funded by the ANR ESPOIR (Edge Simu\-lation of the Physics Of ITER Relevant turbulent transport)and the \emph{F\'ed\'eration nationale de Recherche Fusion par Confinement Magn\'etique} (FR-FCM). We thank Eric Serre, Fr\'ed\'eric Schwander, Guillaume Chiavassa, Philippe Ghendrih and Patrick Tamain for fruitful discussions.

\bibliographystyle{spmpsci}
\bibliography{fvca7-auphan-template}
\end{document}